\documentclass[11pt]{article}

\usepackage{amssymb,amsmath,amsfonts,framed,xcolor}
\usepackage[colorlinks=true, a4paper=true, pdfstartview=FitV,
linkcolor=blue, citecolor=blue, urlcolor=blue]{hyperref}
\usepackage[pdftex]{graphicx}

\newtheorem{thm}{Theorem}[section]
\newtheorem{dfn}[thm]{Definition}

\newtheorem{cor}[thm]{Corollary}

\definecolor{shadecolor}{gray}{0.95}
\newcommand{\HRule}{\rule{\linewidth}{0.5mm}}

\begin{document}

\title{\bfseries{The multicovering radius problem for some types of discrete structures}}
\author{Alan J. Aw}
\date{}
\maketitle

\HRule
\begin{abstract}

The covering radius problem is a question in coding theory concerned with finding the minimum radius $r$ such that, given a code that is a subset of an underlying metric space, balls of radius $r$ over its code words cover the entire metric space. Klapper (\cite{klapper}) introduced a code parameter, called the multicovering radius, which is a generalization of the covering radius. In this paper, we introduce an analogue of the multicovering radius for permutation codes (cf. \cite{keevash-ku}) and for codes of perfect matchings (cf. \cite{aw-ku}). We apply probabilistic tools to give some lower bounds on the multicovering radii of these codes. In the process of obtaining these results, we also correct an error in the proof of the lower bound of the covering radius that appeared in \cite{keevash-ku}. We conclude with a discussion of the multicovering radius problem in an even more general context, which offers room for further research.
\end{abstract}
\HRule

\textbf{Keywords:} probabilistic methods, coding theory, Lov\'asz local lemma, multicovering radius\bigskip

\textbf{Mathematics Subject Classification (2010):} 05D40, 94B99

\section{Introduction}

In this paper, the words \textit{family} and \textit{set} are used interchangeably to denote a collection of sets. We are concerned only with two types of discrete structures: permutations and perfect matchings on simple, finite graphs with no loops. We use the term $m$-set to describe a generic set containing $m$ elements. Let $S_n$ denote the symmetric group acting on the set $[n]$ of positive integers from $1$ to $n$. We use standard graph-theoretic notation, such as those adopted by Bollob\'as in \cite{bollobas}. In particular, let $K_{2n}$ be the complete graph with $2n$ vertices; that is, every two vertices are connected by an edge. A perfect matching, denoted by $M$ throughout, is a graph in which every vertex is incident to exactly one edge. A perfect matching can be represented by a collection of two-element sets where the elements of each set are two distinct vertices. For instance $\{\{v_1,v_2\},\{v_3,v_4\}\}$ is a perfect matching contained in $K_4$.\\

Consider the following problem. Given a metric space $(\Omega,d)$ and any subset $S\subseteq \Omega$ (we call $S$ a \textit{code} in $\Omega$), what is the minimum radius $r$ such that balls of radius $r$ over the points (i.e., code words) in $S$ cover $\Omega$ entirely?\\

This problem is also known as the \textit{covering radius problem}. It is fundamentally related to problems in coding and information theory (see \cite{hill,cohen-honkala} for a comprehensive coverage of the subject, and \cite{cohen} for a survey of recent results). Recently, Cameron and Wanless \cite{cameron} considered the covering radius for subsets $G$ of the symmetric group $S_n$. (They have described the problem to be motivated by a question due to K\'ezdy-Snevily.) In the same vein, the author, together with his research supervisor, \cite{aw-ku} considered the covering radius for subcollections $\mathcal M$ of $1$-factors of a complete uniform hypergraph on $tn$ vertices. A survey by Quistorff \cite{quistorff} provides a concise history and summary of a series of covering and packing problems, of which the covering radius problem is an important aspect.\\ 

In 2006, Ku and Keevash \cite{keevash-ku} introduced a probabilistic approach to establish a frequency parametric result for the covering radius problem for sets of permutations which in some instances yielded tighter lower bounds on covering radius values as compared to the bounds obtained by Cameron and Wanless in \cite{cameron}, thus providing a reliable alternative method to traditional algebraic tools towards investigating the covering radius problem.\\ 

Despite the purely mathematical motivations for these studies, it turns out that the results have potential applications in coding theory or communication models, whereby the applicability depends upon the existence of codes whose behaviour (over a channel) can be characterized by specific discrete structures. For instance, power line communications, a communication model whose data is characterized by permutation arrays \cite{colbourn}, could see applications of results vis-\`a-vis the covering radius problem for sets of permutations. The applications to coding theory notwithstanding, the authors of \cite{keevash-ku} have also applied their results to Latin squares and Latin transversals as it turns out that for certain sets of permutations, the determination of the covering radius is closely related to two important conjectures on Latin squares: Ryser's conjecture that every Latin square of odd order has a transversal, and Brualdi's conjecture that every Latin square of order $n$ has a partial transversal of size $n-1$. Recent literature suggests several possible generalizations of this theory to general groups too.\\   

In this paper, we generalize the covering radius problem for both sets of permutations and sets of perfect matchings to their respective ``\textit{multicovering radius} type" problems. The multicovering radius problem is a generalization of its covering radius counterpart and was, to the best of the author's knowledge, first introduced and studied in significant detail by Klapper\footnote{Summary of the multicovering radius problem on Klapper's homepage: \url{http://www.cs.uky.edu/~klapper/multicov.html}} (see \cite{klapper}). The investigation was motivated by sequence attacks in error correcting codes. Here, an error correcting code of length $n$ is a set $C$ of vectors in the $n$-dimensional space over $GF(2)$. The covering radius of a code $C$ is the smallest integer $r$ such that every vector in $GF(2)^n$ is within distance $r$ of at least one code word in $C$. Here, the distance $d(c_1,c_2)$ between two code words $c_1$ and $c_2$ is just the number of coordinates in which they differ, e.g., $d(101,100)=1$ in $GF(2)^3$. Klapper introduces the general multicovering radius as follows. Let $m\in\mathbb N$. The $m$-covering radius (hence \textit{multi}covering radius) of a code $C$ is the smallest integer $r$ such that for every $m$-set $\{v_1,v_2,...,v_m\}$ of vectors in $GF(2)^n$, there is a vector $c$ in $C$ such that the distance from every $v_i$ to $c$ is at most $r$.\\            

Observing how the covering radius, a notion originating from coding and information theory, became generalized by combinatorialists in \cite{cameron,quistorff,keevash-ku}, it is indeed interesting to ask if it is also possible to generalize the multicovering radius in a similar fashion. Moreover, in view that probabilistic methods have been able to yield stronger results than classical covering bounds, e.g., the sphere-packing bound \cite{quistorff,cameron}, it becomes natural for the author to apply the same tools to generalize the original results (for both sets of permutations and sets of perfect matchings) to obtain their multivariate versions in relation to the multicovering radius problem.\\

On the other hand, while working on this paper we spotted an error \cite{personal} in the proof of the covering radius result in \cite{keevash-ku}. We fix this error by introducing a modified injection, which we also apply to prove our multicovering radius result.\\  

The remaining sections of this paper are divided as follows. First, we introduce the probabilistic tools used to establish our results. Next, we establish results for the multicovering radius problem for collections of permutations which we show are consistent with those obtained for the covering radius problem by Keevash and Ku in \cite{keevash-ku}. Then, we follow the same procedure in the study of the problem for collections of perfect matchings. Finally, we highlight some limitations of the results and briefly discuss some possible analogues of both the covering radius and multicovering radius problems for other types of mathematical structures, which provides room for further research. Remark that we describe the error in \cite{keevash-ku}, as well as its correction, in the appendix.     

\section{Probabilistic Tools}

The main theorem which we use in the paper is the Lov\'{a}sz local lemma. However, for illustrative purposes, elementary results concerning the multicovering radius problem are provided for each type of discrete structure -- permutations and perfect matchings -- studied. These simple results utilize a basic probabilistic result known as the union bound, which we describe below.

\begin{thm}[Union Bound]
\label{union}

Let $\mathcal F$ be a Borel field on a sample space $\Omega$. Then, for a sequence of events $\{A_i\}_{i=1}^{n}, A_i\in\mathcal F$,

\begin{equation*}
\mathbb P\left(\bigcup_{i=1}^{n}{A_i}\right)\leq\sum_{i=1}^{n}{\mathbb P(A_i)}.
\end{equation*}
\end{thm}

Note that the union bound implies the following corollary which will be directly applied to our problems later to yield the elementary results.

\begin{cor}
\label{union2}
Let $A_1,A_2,...,A_n$ be events in a probability space $(\Omega,\mathcal F,\mathbb P)$. If $\sum_{i=1}^{n}{\mathbb P(A_i)}<1$, then

\begin{equation*}
\mathbb P\left(\bigcap_{i=1}^{n}{\overline{A_i}}\right)>0.
\end{equation*} 
\end{cor}

The Lov\'{a}sz local lemma is a powerful tool for showing the existence of structures with desired properties. Briefly speaking, we toss our events onto a probability space and evaluate the conditional probabilities of certain bad events occurring. If these probabilities are not too large in value, then with positive probability none of the bad events occur. More precisely,

\begin{thm}[Lov\'{a}sz; cf. \cite{keevash-ku}, Section 2]
\label{Lovasz}
Let $\mathcal A=\{A_1,\ldots,A_n\}$ be a collection of events in an arbitrary probability space $(\Omega,\mathcal F,\mathbb P)$. A graph $G=(V,E)$ on the set of vertices $V=[n]$ is called a dependency graph for the events $A_1,\ldots,A_n$ if for any $i,j\in V$, $e_{ij}\in E$ iff $A_i$ and $A_j$ are related by some property $\pi$. Suppose that $G$ is a dependency graph for the above events and let $x_1,\ldots,x_n\in[0,1)$. Moreover, for any fixed $i\in [n]$ let $S$ be a subset of the set $[n]\setminus\{j:e_{ij}\in E\}$.  If for each $i$ (and choice of $S$),

\begin{equation*}
\mathbb P\left(A_i\mid\bigcap_{k\in S}{\overline{A_k}}\right)\leq x_i\cdot\prod_{e_{ij}\in E}{(1-x_j)},
\end{equation*}

\noindent then $\mathbb P(\bigcap_{i=1}^{n}{\overline{A_i}})\geq\prod_{i=1}^{n}{(1-x_i)}$. In particular, with positive probability none of the $A_i$ occurs.
\end{thm}

For the proof of Theorem \ref{Lovasz}, we direct the reader to chapter 5 of \cite{alon} or chapter 19 of \cite{jukna}. In this paper, we use the following special case (cf. \cite{alon}) of Theorem \ref{Lovasz} in our result.

\begin{cor}
\label{Symmetric}
Suppose that $\mathcal A=\{A_1,\ldots,A_n\}$ is a collection of events, and for any $A_i\in\mathcal A$ there is a subset $\mathcal D_{A_i}\subset\mathcal A$ of size at most $d$, such that for any subset $\mathcal S\subset\mathcal A\setminus\mathcal D_{A_i}$ we have $\mathbb P\left(A_i\mid\bigcap_{A_j\in\mathcal S}{\overline{A_j}}\right)\leq p$. If $ep(d+1)\leq 1$, then $\mathbb P(\bigcap_{i=1}^{n}{\overline{A_i}})>0$.
\end{cor}

Here, $e:=\lim_{n\to\infty}{\left(1+\frac{1}{n}\right)^n}$.
 
\section{Multicovering Radius of Sets of Permutations}

Consider the permutation group $S_n$ acting on the set $[n]$ of natural numbers from $1$ to $n$. In any collection $G$ of permutations (not necessarily a subgroup), we can measure the {\it Hamming distance} (or distance as we usually drop the first name for brevity) $d(g,h)$ between a permutation $g$ in $G$ and any permutation $h$ picked from $S_n$. Here, the Hamming distance between two permutations is the number of positions in which they differ. For example, in $S_3$, $d(123,231)=3$. It is easy to verify that the Hamming distance for permutations in $S_n$ creates a metric space, and we call $S_n$ endowed with the Hamming metric the {\it Hamming permutation space}. If we were to fix $h$ above and measure the distances $d(h,p)$ for every $p\in G$, there exists a minimum distance which we can obtain between $h$ and some\footnote{There may exist more than one choice of $p_0$ which gives a minimum distance.} $p_0\in G$, i.e., $\min\{d(h,p):p\in G\}=d(g,p_0)$. Now, repeating this procedure for every permutation $h\in S_n$, we can find the maximum of all the minimum distances measured earlier. This maximum value, denoted $cr(G)$, is the covering radius of the collection $G$. Mathematically, $cr(G):=\max_{h\in S_n}\min_{g\in G}{d(g,h)}$. As discussed earlier, there are practical and theoretical motivations towards studying the covering radius problem for such sets.\\

Let us formally define the multicovering radius problem for sets of permutations. For a given subset $G$ of $S_n$ and any $m$-set $\boldsymbol{\Omega}_m$ of permutations in $S_n$ ($m,n\in\mathbb N,m\leq n$), it is possible to compare the distances between each element in $G$ and each element in $\boldsymbol{\Omega}_m$. In particular, for a fixed $m$-set $\boldsymbol{\Omega}_m$, pick any permutation $g$ from $G$ and compute its distance from each permutation in $\boldsymbol{\Omega}_m$. Among the measured distances, take the maximum of them. Now, repeat this procedure for every permutation in $G$ by varying the choice of $g$, and take the minimum of all the maximum measured distances. Finally, vary $\boldsymbol{\Omega}_m$ and repeat the steps above; this gives us a series of ``minimum of maximum'' distances. Among these, pick the maximum value, denoted $cr_m(G)$. Then, $cr_m(G)$ is called the $m$-covering radius of the subset $G$ of $S_n$. Mathematically,

\begin{equation*}
cr_m(G):=\max_{\boldsymbol{\Omega}_m}\min_{g\in G}\max_{h\in\boldsymbol{\Omega}_m}d(g,h).
\end{equation*}

Notice that when $m=1$, the formula reduces to

\begin{equation*}
cr_1(G)=cr(G)=\max_{h\in S_n}\min_{g\in G}d(g,h).
\end{equation*}       

This is indeed the definition of the covering radius of $G\subseteq S_n$ as described earlier. Moreover, it is clear that such a definition of the multicovering radius is consistent with the one introduced by Klapper as described in the Introduction.

\subsection{Results}

An important question to ask is: given any collection of permutations $G\in S_n$, what is its $m$-covering radius? For small $n$, $cr_m(G)$ can be computed easily. However, for general $n$, it remains an unsolved problem to accurately determine $cr_m(G)$ for any $G\subseteq S_n$. Even for the covering radius problem, i.e., the case $m=1$, there exists no known explicit formula for the covering radius as $n$ gets large. Thus, estimations of the $m$-covering radius are established instead.\\

First, we establish, using corollary \ref{union2}, a lower bound for the $m$-covering radius. This is for illustrative purposes and is a generalization of an elementary result described in \cite{keevash-ku}. This is expected to be weaker than other possibly existent bounds, e.g., the sphere-packing bound established by Cameron and Wanless which can be extended for the $m$-covering radius. We suggest that the interested reader should refer to \cite{cameron,quistorff}, inasmuch as the main purpose here is to demonstrate the use of the probabilistic method. 

\begin{thm}
\label{permutation1}
Let $G\subseteq S_n$ such that $|G|<\frac{{n!\choose m}}{{(n-s)!\choose m}\cdot{n\choose s}}$. Then $cr_m(G)\geq n-s+1$.
\end{thm}

\noindent\textbf{Proof.}\\

Let $G=\{g_1,g_2,...,g_k\}$. Pick an $m$-set $\boldsymbol{\Omega}_m$ uniformly at random. Given an index $i\in\{1,2,...,k\}$ and a set $S\subset [n]$ of size $s$, define $A_{i,S}$ to be the event that all permutations in $\boldsymbol{\Omega}_m$ agree with $g_i$ on $S$, i.e., $\forall g\in\boldsymbol{\Omega}_m,~g(x)=g_i(x)$ for all $x\in S$. Then, $\mathbb P(A_{i,S})=\frac{{(n-s)!\choose m}}{{n!\choose m}}$ since there are exactly $(n-s)!$ permutations which have $s$ common positions. Summing over all $(i,S)$, it follows that

\begin{equation*}
\sum_{(i,S)}{\mathbb P(A_{i,S})}=k{n\choose s}\frac{{(n-s)!\choose m}}{{n!\choose m}}<1.
\end{equation*}

So, by corollary \ref{union2} there exists an $m$-set which agrees with every element of $G$ in at most $s-1$ positions. This implies $cr_m(G)\geq n-s+1$.\hfill $\blacksquare$\\

We now establish a frequency parametric lower bound for the $m$-covering radius in terms of the following frequency parameter: for $G\subseteq S_n$ and $1\leq a,b\leq n$, let $N_G(a,b)=|\{g\in G:g(a)=b\}|$. Notice that the values of $N_G(a,b)$ for pairs $a,b\in [n]$ impose some restrictions on the size of $G$. In this proof, we apply corollary \ref{Symmetric} using a strategy which is similar to that used in the proof of the Erd\H{o}s-Spencer theorem on Latin transversals, as presented in chapter 5 of \cite{alon} (pp. 73-74). We suggest that readers who are new to the Lov\'asz local lemma read through that proof (i.e., Erd\H{o}s-Spencer in \cite{alon}). Note that this type of lower bound was established for the covering radius problem in \cite{keevash-ku} and was shown to yield better computational results, both in terms of computational speed and tightness of bound, than the standard sphere-packing bound for certain values of $n$.  

\begin{thm}
\label{permutation2}

Let $G\subseteq S_n$ be a collection of permutations such that $N_G(a,b)\leq k$ for any $a,b\in [n]$. If

\begin{equation*}
k\leq \frac{(n-s)!}{(n-1)!(2n-s)}\cdot\frac{(s-1)!}{s}\cdot\left[\frac{n!}{(n-s)!}\right]^m\left(\frac{1}{e}-\left[\frac{(n-s)!}{n!}\right]^m\right)
\end{equation*}

\noindent for some positive integer $s$, then there exists an $m$-set of permutations whose elements each agree with each permutation of $G$ in at most $s-1$ positions, i.e., $cr_m(G)\geq n-s+1$. 
\end{thm}   

\noindent\textbf{Proof.}\\

Pick an $m$-set $\boldsymbol{\Omega}_m$ uniformly at random. We shall show that with positive probability each element in $\boldsymbol{\Omega}_m$ agrees with each permutation in $G$ in at most $s-1$ positions.\\

Let $G=\{g_1,g_2,...,g_r\}$. Given an index $i\in\{1,2,...,r\}$ and a set $S\subset [n]$ of size $s$, define $A_{i,S}$ to be the event that all permutations in $\boldsymbol{\Omega}_m$ agree with $g_i$ on $S$. Let $\mathcal A$ be the set of all the events $A_{i,S}$. We also let $X_{i,S}$ be the collection of pairs $(i',S')$ such that at least one of $S\cap S'$ or $g_i(S)\cap g_{i'}(S')$ is non-empty (for a function $f$ and a subset $S$ of its domain, $f(S)=\{f(s):s\in S\}$). Let $\mathcal D_{i,S}$ comprise the events $A_{i',S'}$ such that $(i',S')\in X_{i,S}$. Let us count that number of events $A_{i',S'}\in\mathcal D_{i,S}$. First, choose two elements $x,y\in [n]$ so that at least one of $x\in S$ or $y\in g_i(S)$ holds: there are $2sn-s^2$ choices. Next, choose $i'$ such that $g_{i'}(x)=y$; there are at most $k$ choices by our assumption that $N_G(a,b)\leq k$ for any $a,b\in [n]$. Finally, the rest of $S'$ can be chosen in at most ${n-1\choose s-1}$ ways. Therefore,

\begin{equation*}
|\mathcal D_{i,S}|\leq ks(2n-s){n-1\choose s-1}=d.
\end{equation*}

Let us bound $\mathbb P(A_{i,S}\mid E)$, where $E=\bigcap_{A_{i',S'}\in\mathcal S}\overline{A_{i',S'}}$ for any subset $\mathcal S\subseteq\mathcal A\setminus\mathcal D_{i,S}$. Now, for any $\boldsymbol{\Omega}_m$ picked, it is possible to order its elements lexicographically according to a fixed ordering of the permutations in $S_n$. In what follows, let $f:S\mapsto [n]$ be any injection, and let $B_f$ be the event that a permutation $g\in S_n$ restricts to $f$ on $S$. Let the set of all such $f$ be $F$. Consider any $m$-tuple $\boldsymbol{f}$ contained in $F^m$. Let $B_{\boldsymbol{f}}$ be the event that, for component $i~(1\leq i\leq m)$ (for brevity, call it $f_i$) in $\boldsymbol{f}$, the $i$th permutation (with respect to the lexicographic order imposed) in $\boldsymbol{\Omega}_m$ restricts to $f_i$ on $S$.\\

\noindent\textit{Claim: $\mathbb P(A_{i,S}\mid E)\leq\mathbb P(B_{\boldsymbol{f}}\mid E)$.}\\

The claim is established by means of constructing an injective map from the collection of $m$-sets $\boldsymbol{\Omega}_m$ such that $A_{i,S}\cap E$ holds, to the collection of $m$-sets such that $B_{\boldsymbol{f}}\cap E$ holds. This map is to replace each $\boldsymbol{\Omega}_m$ satisfying $A_{i,S}\cap E$ with some unique $\boldsymbol{\Omega}_m'$ satisfying $B_{\boldsymbol{f}}\cap E$. To do so, we order the permutations in $\boldsymbol{\Omega}_m$ lexicographically as mentioned earlier, and then perform a mapping on every permutation in $\boldsymbol{\Omega}_m$ such that the $i$th permutation is mapped to another permutation $h$ which restricts to $f_i$ on $S$ for $1\leq i\leq m$. The image of such a mapping is another $m$-set $\boldsymbol{\Omega}_m'$, but one which clearly satisfies $B_{\boldsymbol{f}}\cap E$.\\

Indeed, this map is to replace the $i$th permutation in a $\boldsymbol{\Omega}_m$ that satisfies $A_{i,S}\cap E$, which we denote here by $h_i$, to another permutation $h$, which is defined by a mapping $\phi$ described below.\footnote{The mapping here is different from the original mapping used by Keevash and Ku in \cite{keevash-ku} which, we show in the appendix, is problematic.}

\begin{framed}
\noindent\textbf{The $\phi:h_i\mapsto h$ mapping.}\\
\textbf{1.} Let $h(x)=f_i(x)$ for $x\in S$.\\
\textbf{2.} Consider the set $T=\{x\in [n]\setminus S:h_i(x)\in f_i(S)\}$. Define the composite function $\rho=h_i(f_i^{-1}(\cdot))$. For all $x\in T$, let $h(x)=\rho^{(N)}(h_i(x))$, where $N$ is the minimum positive integer such that $\rho^{(N)}(h_i(x))\not\in f_i(S)$.\\
\textbf{3.} Let $h(x)=h_i(x)$ for all other $x\not\in S\cup T$.
\end{framed}

First, it is not too difficult (albeit not immediate) to show\footnote{The appendix provides a detailed proof.} that $h(x)\neq h(y)$ as long as $x\neq y$, by applying the fact that $h_i$ and $f_i$ are bijective. This ensures that $h$ is indeed a permutation. Second, observe that $\phi$ is injective. Indeed, suppose $h_i$ and $h_i'$ are distinct and come from the same components of two distinct lexicographically ordered $\boldsymbol{\Omega}_m$ that satisfy the events $A_{i,S}$ and $E$. Thus there is some position $x_0\not\in S$ in which $h_i$ and $h_i'$ differ, since by hypothesis $h_i(S)=h_i'(S)$. Suppose $x_0\in T$, or else we are done. Since $x_0\in T$, observe that $\rho$ is a composition of two bijective mappings, and that $h_i(x)=h_i'(x)~~\forall x\in S$; this ensures that $h(x_0)\neq h'(x_0)$. (It does not matter even if $M$ and $N$ are different.) Hence, upon repeating $\phi$ for every $h_i\in\boldsymbol{\Omega}_m$, an $m$-set is obtained which can be verified to satisfy both $B_{\boldsymbol{f}}$ and $E$. (Our definition of $X_{i,S}$ ensures that none of the events in $E$ is affected by $\phi$, since $\phi$ only affects $S$ and $h_i(S)$.)\\

Therefore, inasmuch as the events $B_{\boldsymbol{f}}\cap E$ are mutually exclusive for different $\boldsymbol{f}$, and $\mathbb P\left(\bigcup_{\boldsymbol{f}}{B_{\boldsymbol{f}}\cap E}\right)=\mathbb P(E)$, we have

\begin{equation*}
\mathbb P(A_{i,S}\mid E)\leq \left[\frac{(n-s)!}{n!}\right]^m\sum_{\boldsymbol{f}}{\mathbb P(B_{\boldsymbol{f}}\mid E)}=\left[\frac{(n-s)!}{n!}\right]^m=p.
\end{equation*}            

Now, we want $ep(d+1)\leq 1$. This is equivalent to our bound on $k$.\hfill $\blacksquare$\\

It is worth mentioning that when $m=1$, the result above is exactly the one obtained in \cite{keevash-ku}.

\section{Multicovering Radius of Sets of Perfect Matchings} 

The covering radius problem for sets of perfect matchings was first studied in \cite{aw-ku}. However, it is not the first instance in which the covering radius problem was studied with respect to graphs. Quistorff in \cite{quistorff} highlighted the study of this problem in relation to distances between two vertices of a graph. In graph theoretic terminology, it is known as the $e$-domination number. The interested reader should read \cite{quistorff} and relevant literature regarding that problem.\\

For our problem, which is different from the one just described above, we are considering collections $\mathcal M$ of perfect matchings in the complete graph $K_{2n}$ on an even number of vertices. We work in the finite metric space $(\Omega,d)$ where $\Omega$ is the set of all perfect matchings of $K_{2n}$ and $d$ is defined for any two matchings $M,M'$ to be the number of edges in which they differ, i.e., 

\begin{equation*}
d(M,M')=n-|M\cap M'|
\end{equation*}

\noindent For example, in $K_6$, the perfect matchings $M_1=\{\{v_1,v_2\},\{v_3,v_4\},\{v_5,v_6\}\}$ and $M_2=\{\{v_1,v_3\},\{v_2,v_4\}$ 

,$\{v_5,v_6\}\}$ satisfy $d(M_1,M_2)=2$. Indeed, we could think of $d(M,M')$ as the ``Hamming distance" between $M$ and $M'$. Now, given a collection $\mathcal M$ of perfect matchings, its covering radius is defined as follows: fix a perfect matching $M'$ of the universal set of all perfect matchings in $K_{2n}$, and measure the distances $d(M,M')$ for every $M\in\mathcal M$. Pick the minimum distance out of all distances measured. Repeat this procedure for each perfect matching of the set of all perfect matchings in $K_{2n}$, and a series of minimum distances is obtained. Then, the maximum value is the covering radius of $\mathcal M$, denoted $cr(\mathcal M)$. Thus, the covering radius problem could be stated as follows: given a collection $\mathcal M$ of perfect matchings in $K_{2n}$, what is the largest possible number of elements in this collection such that we can find a perfect matching in $K_{2n}$ that agrees with each perfect matching in the collection in at most $x-1$ edges?\\

The notion of $m$-covering radii of sets of perfect matchings is similar to that of the $m$-covering radii of sets of permutations. Indeed, given a collection $\mathcal M$ of perfect matchings, pick one of its elements, say $M$. Consider all $m$-sets $\boldsymbol{\Lambda}_m$ of perfect matchings, i.e., sets containing $m$ perfect matchings; and select one particular $m$-set, $\boldsymbol{\Lambda}_m$. Compute the distance $d(M,M')$ between $M\in\mathcal M$ and each $M'$ from $\boldsymbol{\Lambda}_m$. This gives us a series of distances, of which we keep the maximum. Now, vary $M$ to obtain a sequence of maximum distances. Take the minimum of these distances. Lastly, vary the choice of $\boldsymbol{\Lambda}_m$, and repeat the procedure outlined above. Upon obtaining a series of ``minimum of maximum'' distances, pick the maximum; this value is the $m$-covering radius of $\mathcal M$, denoted $cr_m(\mathcal M)$. Indeed, mathematically

\begin{equation*}
cr_m(\mathcal M)=\max_{\boldsymbol{\Lambda}_m}\min_{M\in\mathcal M}\max_{M'\in\boldsymbol{\Lambda}_m}d(M,M').
\end{equation*}

It can be checked easily that when $m=1$, the formula reduces to one which is consistent with our definition of the covering radius of a set of perfect matchings in $K_{2n}$.

\subsection{Results}

Here, we introduce a term called $x$-matching.

\begin{dfn}
\label{matching}
An $x$-matching of $K_{2n}$ is a matching of size $x$. Thus, if $x=n$, then the $x$-matching is simply a perfect matching.
\end{dfn}

Additionally, for $x\leq y$, we say that a $y$-matching $Y$ \textit{contains} a $x$-matching $X$ if and only if all the edges in $X$ are also edges in $Y$ ($X\subseteq Y$). We begin by providing a basic (probabilistic) result, again for illustrative purposes. From corollary \ref{union2}, we obtain an elementary bound concerning the multicovering radius of a set of perfect matchings.

\begin{thm}
\label{pm1}

Let $\mathcal M$ be a collection of perfect matchings in $K_{2n}$. Moreover, denote $\beta_x=\frac{(2n-2x)!}{2^{n-x}\cdot(n-x)!}$ where $1\leq x\leq n$; additionally, let $\beta=\frac{(2n)!}{2^{n}\cdot n!}$. If $|\mathcal M|<\frac{{\beta\choose m}}{{\beta_x\choose m}\cdot{n\choose x}}$, then $cr_m(\mathcal M)\geq n-x+1$.
\end{thm}

\noindent\textbf{Proof.}\\

Let $\mathcal M=\{M_1,...,M_s\}$. Pick an $m$-set $\boldsymbol\Lambda_m$ independently, uniformly and randomly. Given an index $i\in\{1,2,...,m\}$ and a $x$-matching $X\subset M_i$, let $A_{i,X}$ be the event that all perfect matchings in $\boldsymbol\Lambda_m$ contain $X$, i.e., $\forall M\in\boldsymbol\Lambda_m, X\in M$. Moreover, note that there are $\frac{(2n)!}{2^n\cdot n!}=\beta$ perfect matchings in $K_{2n}$, and exactly $\frac{(2n-2x)!}{2^{n-x}\cdot(n-x)!}=\beta_x$ of them with $x$ fixed edges. Then, $\mathbb P(A_{i,X})=\frac{{\beta_x\choose m}}{{\beta\choose m}}$. Summing over all $(i,X)$, inasmuch as there are $s$ possible values of $i$, and for each perfect matching $M_i$ there are exactly ${n\choose x}$ $x$-matchings, it follows that 

\begin{equation*}
\sum_{(i,X)}{\mathbb P(A_{i,X})}=s{n\choose x}\frac{{\beta_x\choose m}}{{\beta\choose m}}<1.
\end{equation*}

So, by corollary \ref{union2} there exists an $m$-set which agrees with every perfect matching of $\mathcal M$ in at most $x-1$ positions. This implies that $cr_m(\mathcal M)\geq n-x+1$.\hfill $\blacksquare$\\ 

Next, we establish a frequency parametric result, which is similar to Theorem \ref{permutation2}. This time, the frequency parameter counts the number of times each edge appears in an explicit listing of elements $M\in\mathcal M$. Again, the proof is an application of corollary \ref{Symmetric}, and the strategy used is similar to the one used to establish Theorem \ref{permutation2}.

\begin{thm}
\label{pm2}

Let $\mathcal M$ be a collection of perfect matchings in $K_{2n}$ such that each of the ${2n\choose 2}$ edges appears at most $k$ times in an explicit listing of the elements $M_i\in\mathcal M$. If

\begin{equation*}
k\leq\frac{1}{2x(2n-1){n-1\choose x-1}}\cdot\left[\sum_{k=2x-n}^{x}{{2x\choose 2k}\frac{(2k)!}{k!\cdot2^k}}\right]^m\cdot\left(\frac{1}{e}-\left[\sum_{k=2x-n}^{x}{{2x\choose 2k}\frac{(2k)!}{k!\cdot2^k}}\right]^{-m}\right)
\end{equation*}

\noindent for some positive integer $x$, then there exists an $m$-set of perfect matchings whose elements each agree with each perfect matching $M_i\in\mathcal M$ in at most $x-1$ edges, i.e., $cr_m(\mathcal M)\geq n-x+1$.
\end{thm}

\noindent\textbf{Proof.}\\

Pick an $m$-set $\boldsymbol\Lambda_m$ randomly, independently, and uniformly. We shall show that with positive probability each element in $\boldsymbol\Lambda_m$ agrees with each perfect matching in $\mathcal M$ in at most $x-1$ positions.\\

Let $\mathcal M=\{M_1,...,M_s\}$. Given an index $i\in\{1,2,...,s\}$ and a $x$-matching $X\subset M_i$, define $A_{i,X}$ to be the event that all perfect matchings in $\boldsymbol\Lambda_m$ contain $X$. Let $\mathcal A$ be the set of all the events $A_{i,X}$. We also let $Q_{i,X}$ be the collection of pairs $(i',X')$ such that $X$ and $X'$ share at least one common vertex in their underlying vertex sets. Let $\mathcal D_{i,X}$ comprise the events $A_{i',X'}$ such that $(i',X')\in Q_{i,X}$. Let us count the number of events $A_{i',X'}\in\mathcal D_{i,X}$. First, pick one vertex out of the $2x$ vertices in $X$, then choose out of the $2n-1$ remaining vertices of $M_i$ one particular vertex to be its neighbour in the perfect matching. Next, choose a perfect matching $M_{i'}\in \mathcal M$ that contains the constructed edge; this can be done in at most $k$ ways by our assumption. Finally, the rest of $X'$ can be chosen in ${n-1\choose x-1}$ ways. Therefore

\begin{equation*}
|\mathcal D_{i,X}|\leq k\cdot2x(2n-1){n-1\choose x-1}=d.
\end{equation*}

Let us now consider the probability $\mathbb P(A_{i,X}\mid E)=p_0$, where $E=\bigcap_{A_{i',X'}\in\mathcal S}\overline{A_{i',X'}}$ for any subset $\mathcal S\subseteq\mathcal A\setminus\mathcal D_{i,X}$. For the rest of this proof, our aim is to bound $p_0$ from above.\\

Now, for any $\boldsymbol\Lambda_m$ picked, it is possible to order its elements lexicographically, according to a fixed ordering of the perfect matchings in $K_{2n}$ that gives us a bijection between the set of natural numbers from $1$ to $\frac{(2n)!}{n!\cdot 2^n}$ and $\Omega$. We exploit this fact and treat $\boldsymbol\Lambda_m$ as a ``lexicographically ordered'' vector with $m$ components where each component is a perfect matching. In what follows, we shall describe a procedure to bound the set of $\boldsymbol\Lambda_m$ satisfying $A_{i,X}\cap E$.\\

Fix $A_{i,X}$. Without loss of generality, let the underlying set of $2x$ vertices of $X$ be $V=\{v_1,v_2,...,v_{2x}\}$ and the rest of the vertices not in $X$ be $\{v_{2x+1},..., v_{2n}\}$. Arbitrarily partition the $2x$ vertices contained in $X$ into a collection $W$ of edges and vertices, additionally treating a vertex as a vertex contained in a set, i.e., a singleton. By this, we mean, for example, that if $X=\{\{v_1,v_2\},\{v_3,v_4\}\}$ then a possible partition $W$ is $\{\{v_1\},\{v_2\},\{v_3,v_4\}\}$ instead of $\{v_1,v_2,\{v_3,v_4\}\}$. Letting the number of singletons (sets each containing a vertex) be $2q$ and the number of doubletons (edges) be $r$ (by virtue that the union of the partition has an even-numbered size), we further restrict our partitions to only partitions such that $2q+r\leq n$. Note that $q+r=x$. The reason we make such a restriction is because it gives us $n-x \geq p \geq0$. Its importance will become apparent later.\\

Let $\mathcal W$ be the family of all such $W$. Consider any $m$-tuple $\boldsymbol{W}\in\mathcal W^m$. Let $B_{\boldsymbol{W}}$ be the event that, for component $i$ - where $1\leq i\leq m$ - (for brevity, call it $W_i$) in $\boldsymbol{W}$, the $i$th perfect matching (with respect to the lexicographic order imposed) in $\boldsymbol\Lambda_m$ \textit{contains} $W_i$. By containment, we mean, for example, that $\{\{v_1,v_2\},\{v_3,v_4\},\{v_5,v_6\}\}$ contains $\{\{v_1\},\{v_3\},\{v_5,v_6\}\}$.\\

\noindent\textit{Claim: $\mathbb P(A_{i,X}\mid E)\leq\mathbb P(B_{\boldsymbol{W}}\mid E)$}.\\

The claim is established by means of constructing an injective map from the collection of $m$-sets $\boldsymbol\Lambda_m$ such that $A_{i,X}\cap E$ holds, to the collection of $m$-sets such that $B_{\boldsymbol{W}}\cap E$ holds. This map is to replace each $\boldsymbol\Lambda_m$ satisfying $A_{i,X}\cap E$ with some unique $m$-set satisfying $B_{\boldsymbol{W}}\cap E$. To do so, we treat every $\boldsymbol\Lambda_m$ as a vector, and then perform a mapping on each component such that the $i$th component is mapped to another perfect matching which contains $W_i$ for $1\leq i\leq m$. The image of such a mapping is another vector which, upon removal of the imposed lexicographic order, is effectively another $m$-set, but one which satisfies $B_{\boldsymbol{W}}\cap E$.\\

Let us denote our mapping on the vector $\boldsymbol\Lambda_m$ (satisfying $A_{i,X}\cap E$) by $\phi$. Since $\phi$ is to be performed identically on each component of the vector, we just describe $\phi$ in detail for the first component, which we denote by $M_1$ throughout our explanation for brevity. First, before defining $\phi$ explicitly, consider the $n-x$ edges of $M_1$ not in $X$. Order the vertices in each of these edges based on the natural ordering of the vertices' subscripts. Thus, every edge not in $X$ now becomes a directed edge $(v_{i_1},v_{i_2})$ where $i_1<i_2$. Next, order these $n-x$ directed edges by the following rule: compare every two directed $t$-edges and place the one whose first element is a vertex with a smaller subscript on the left hand side of the other; this gives a sequence of $n-x$ directed edges

\begin{equation*}
(v_{\tau_{M_1}(1)},v_{\tau_{M_1}(2)}),(v_{\tau_{M_1}(3)},v_{\tau_{M_1}(4)}),\ldots,(v_{\tau_{M_1}(2n-2x-1)},v_{\tau_{M_1}(2n-2x)}),
\end{equation*} 

\noindent where $\tau_{M_1}$ is a bijection from $\{1,2,\ldots,2n-2x\}$ to $\{2x+1,\ldots,2n\}$ such that

\begin{equation*}
\tau_{M_1}(1)< \tau_{M_1}(3)<\cdots<\tau_{M_1}(2n-2x-1).
\end{equation*}

Now consider a fixed $\boldsymbol{W}\in\mathcal W^m$ and its first component $W_1$. Recall that we want to map $M_1$, which contains $X$, to $\phi(M_1)$, which contains $W_1$. Let $W_1$ contain $2p$ singletons and $q$ doubletons. Construct a bijection $\sigma$ between the singletons and $[2p]$ as follows. Order the singletons into a string of length $2p$ such that between any two vertices the one with a smaller subscript appears first. Then, for each $i$, map the $i$th position of the string to the natural number $i\in[2p]$. This yields

\begin{equation*}
W^{(1)}_1=\{v_{\sigma(1)}\},W^{(2)}_1=\{v_{\sigma(2)}\},\ldots, W^{(2p)}_1=\{v_{\sigma(2p)}\},
\end{equation*}
      
\noindent where $\sigma(i)<\sigma(j)$ whenever $i<j$, and $\bigcup_{i\in[2p]}\{v_{\sigma(i)}\}$ is just the set of all vertices which are singletons in the partitioning of $X$ into $W_1$.\\

Now construct the following $2p$ edges:

\begin{eqnarray}
E_{1} & = & W^{(1)}_1\cup\{v_{\tau_{M_1}(1)}\},\nonumber\\
	& \vdots & \nonumber\\
E_{i} & = & W^{(i)}_1\cup\{v_{\tau_{M_1}(i)}\},\nonumber\\
	& \vdots & \nonumber\\
E_{2p} & = & W^{(2p)}_1\cup\{v_{\tau_{M_1}(2p)}\}.\nonumber
\end{eqnarray}   

Note that $2p\leq2n-2x$ so the vertices $v_{\tau_{M_1}(1)},\ldots,v_{\tau_{M_1}(2p)}$ exist and were used in the above construction. Moreover, inasmuch as an even number of vertices from the set of edges not contained in $X$ were used, there remains an even number of vertices not in $X$ which remain unchanged. These vertices are moreover ordered in relation to their edges such that we can immediately retrieve the respective edges that they were each a part of. We call the family of these retrieved edges $F$. (Note that $F$ is empty if $2p=2n-2x$.)\\

We shall now define our injection $\phi:M_1\mapsto \phi(M_1)$. Set $\phi(M_1)$ to be the perfect matching containing the edges $E_{1},\ldots,E_{2p}$, the doubletons (edges) in $W_1$ that were untouched, and the edges in $F$. Observe that after applying $\phi$ to every component in the vector $\boldsymbol\Lambda_m$, its image, which we shall conveniently represent by $\phi(\boldsymbol\Lambda_m)$, satisfies $B_{\boldsymbol{W}}$. Moreover, since for each component, in particular the first component $M_1$, and every $X¡¯\cap X=\emptyset$, $X¡¯\not\subseteq M_1$ implies $X¡¯\not\subseteq \phi(M_1)$, we conclude that $\phi(\boldsymbol\Lambda_m)$ satisfies $E$.\\    

It remains to show that $\phi$ is injective. We show this by showing that it is injective for each component, particularly the first component. Consider two distinct $m$-sets $\boldsymbol\Lambda_m$ and $\boldsymbol\Lambda'_m$ satisfying $A_{i,X}\cap E$ such that without loss of generality their first components differ, i.e., $M_1\neq M'_1$. Since they agree on $X$, they differ outside $X$ and thus there exists a first edge not in $X$ in which they differ. Suppose this edge is $\{v_{a},v_{b}\}$ for $M_1$ and $\{v_{c},v_{d}\}$ for $M'_1$. Assuming $a<b$ and $c<d$, it is not difficult to observe that either $v_a\neq v_c$ or $v_b\neq v_d$ (or both). After ordering the vertices to yield the directed edges $(v_a,v_b)$ and $(v_c,v_d)$, and applying the bijections $\tau_{M_1}$ and $\tau_{M'_1}$ on the respective edge sets, it is clear that there exists $n\in\mathbb N$ such that $\tau_{M_1}(k)=a,\tau_{M'_1}(k)=c$ and $\tau_{M_1}(k+1)=b,\tau_{M'_1}(k)=d$. Now, if $k>2p$, then we are done because this would mean that the four vertices are not mapped by $\phi$, so that by retrieving these four unmapped vertices back into their original edges, we are guaranteed that $\phi(M_1)$ contains $\{v_a,v_b\}$ while $\phi(M'_1)$ contains $\{v_c,v_d\}$, i.e., $\phi(M_1)\neq\phi(M'_1)$. On the other hand, if $k\leq 2p$, then the edges $E_k=W^{(k)}_1\cup\{v_a\}$ and $E_{k+1}=W^{(k+1)}_1\cup\{v_b\}$ are contained in $\phi(M_1)$ while the edges $E'_k=W^{(k)}_1\cup\{v_c\}$ and $E'_{k+1}=W^{(k+1)}_1\cup\{v_d\}$ are contained in $\phi(M'_1)$. But clearly at least one of the inequalities $E'_k\neq E_k, E'_{k+1}\neq E_{k+1}$ holds. Therefore $\phi(M_1)\neq\phi(M'_1)$ as well, and we conclude that $\phi$ is indeed injective.\\

Hence, we have

\begin{equation*}
p_0=\mathbb P(A_{i,X}\mid E)\leq\mathbb P(B_{\boldsymbol{W}}\mid E).
\end{equation*}

Summing all possible $\boldsymbol W$, and observing that the events $B_{\boldsymbol{W}}$ are mutually exclusive with $\mathbb P\left(\bigcup_{\boldsymbol{W}}{B_{\boldsymbol{W}}}\right)=1$, we have $p_0\cdot N\leq 1$ where

\begin{equation*}
N=\left[\sum_{r=2x-n}^{x}{{2x\choose 2r}\frac{(2r)!}{r!\cdot 2^r}}\right]^m
\end{equation*}

\noindent is the number of $m$-tuples contained in $|\mathcal W^m|$. Notice that each summand counts the number of ways to partition the underlying vertex set into a collection of $r$ doubletons and $2x-2r$ singletons; and we set ${n\choose r}=0$ whenever $r<0$.\\

This gives

\begin{equation*}
p_0\leq \frac{1}{N}=p.
\end{equation*}

In view of corollary \ref{Symmetric}, we want $ep(d+1)\leq1$. This is equivalent to our bound on $k$.\hfill $\blacksquare$\\

It is worth mentioning that when $m=1$, the result above is exactly the one obtained in \cite{aw-ku}.

\section{Concluding Remarks}

In our study of the multicovering radius problem, we depended only on probabilistic tools to derive our results. There may be other methods not known to the author which could yield comparable or even better results. Several potentially potent techniques may be found in \cite[Chapters 30--34]{graham}. Moreover, as discussed in \cite{keevash-ku,aw-ku}, practical applications (both in mathematics and in engineering) of the results could also be considered when evaluating the goodness of these established results. For one thing, we were cavalier in our bounding of the degree of each event in the dependency graph, and this could greatly worsen the calculated results during practical applications. In the event that this happens, one could apply elementary combinatorial methods to improve the crude bounds obtained above. Inasmuch as actual applications sometimes involve data that can provide more than just the frequency parameter $k$ (e.g., the actual size of the collection), combinatorial methods that improve the bound described above will greatly improve the optimality of the output values of covering radius. Apart from these limitations, other questions, particularly the asymptotic behaviour of the formulae obtained, may also be of interest.\\

In this paper, we investigated only discrete structures under the Hamming metric. From our investigation, it is clear that the multicovering radius problem can be seen generally as the minimum radius $r$ such that, given a code $S\subseteq\Omega$, balls of radius $r$ over its code words cover the family of all $m$-sets of $\Omega$. Thus, it is possible to consider other types of metric spaces $(\Omega,d)$ and their possible multicovering and covering analogues (certainly, the notion of distance $d$ is closely related to the class of objects in $\Omega$). For instance, one could consider the multicovering radius for subsets of $C_n$, the cyclic group, or even sets of integers, lattices and partially ordered sets. Sometimes, the metric chosen to act on the class of objects might not be obvious or have any clear motivations. In other instances, e.g., in the latter case, the metric could perhaps be the classical Euclidean distance $|m-n|$, or even the Pythagorean distance $\sqrt{|m^2-n^2|}$. Then, the covering radius could first be defined in an analogous manner, from which it could also be generalized to the multicovering radius in a sensible manner. (For instance, if $[n]$ is the underlying set and $S\subseteq [n]$, then $cr(S)=\max_{b\in [n]}\min_{a\in S} d(a,b)$.) Based on the author's knowledge, there is no literature that has studied the multicovering problem other than for vectors in $\mathbf F^n$ and specifically rank distance (RD) codes \cite{vasantha}.\bigskip
 
\noindent\textbf{Acknowledgments.} The author would like to thank Dr Cheng Yeaw Ku for his guidance and generous help towards the creation of this paper; Prof Peter Keevash and Prof Peter Cameron for their kind advice with regards to improvements and submission of the paper; and the anonymous referee for his helpful comments.

\noindent\textit{Raffles Science Institute\\
Raffles Institution,
One Raffles Institution Lane, S575954,\\
Singapore\\
nalawanij@gmail.com} 

\newpage

\section*{Appendix}

Towards the end of our proof of Theorem \ref{permutation2}, we introduced the mapping $\phi:h_i\mapsto h$, which we claim to be injective. Originally, we intended to use the mapping by Keevash and Ku in \cite{keevash-ku} in establishing our result. Their mapping, with some symbols and notation modified to suit the context of our proof, is as follows.

\begin{framed}
``For any injection $f_i:S\rightarrow [n]$ let $B_{f_i}$ be the event that $h_i$ restricts to $f_i$ on $S$. (...) To see this, we exhibit an injective map, from the set of permutations $h_i$ such that $E\cap A_{i,S}$ holds, to the set of permutations such that $E\cap B_{f_i}$ holds. This map is to replace the permutation $h_i$ by $h$, which is defined as follows. Let $T=\{x\in [n]\setminus S:h_i(x)\in f_i(S)\}$. Define $h(x)=f_i(x)$ for $x\in S$, $h(x)=h_i(f_i^{-1}(h_i(x)))$ for $x\in T$, and $h(x)=h_i(x)$ otherwise. It it not hard to see that the map $h_i\mapsto h$ is injective, and that $E\cap B_{f_i}$ holds for $h$ (using the definition of $X_{i,S}$ to see that none of the events in $E$ is affected by the map)."
\end{framed}

Here, $A_{i,S}$ and $E$ can be seen to be defined in the same way as how we have defined them in our proof; in particular $A_{i,S}$ is just the event that a randomly picked permutation agrees with $g_i\in G$ on $S$ (since the authors of \cite{keevash-ku} were considering it in the context of the covering radius problem and not the multicovering radius problem).\\

In fact, it is not difficult to prove that the map is injective and preserves $E$. The problem lies in that it has not been shown that $h$ is a permutation. Here, we construct an example that shows that $h$ is \textit{not always} a permutation under the mapping described above.\\

Consider $S_7$ and make the following definitions.

\begin{eqnarray}
f_i(123) & = & 123 \\
h_i & = & 7123456 
\end{eqnarray}

Based on these definitions, $S=\{1,2,3\}$. Under the mapping described by the authors, we yield $h=1232456$, where $h(4)=h_i(f_i^{-1}(h_i(4)))=g(3)=2$. Clearly $h$ is not a permutation.\\

Despite so, the mapping described above can be modified to yield a $h$ that is always a permutation, while also fulfilling the rest of the assumptions. This mapping is the one described in our proof. Below, we restate our mapping and provide a detailed argument explaining why it works.

\subsection*{Fixing the injective map}

The frame below provides a detailed procedure of our mapping, which is slightly modified from the mapping by the authors of \cite{keevash-ku}. Here, we stick to the definitions of $A_{i,S}$ and $E$ in the appendix to keep the flow of the presentation and argument. Note that it should not be difficult at all to absorb the meaning of these definitions because they are almost identical to the ones defined in our proof of Theorem \ref{permutation2}.  

\begin{framed}
\noindent\textbf{The $\phi:h_i\mapsto h$ mapping.}\\
\textbf{1.} Let $h(x)=f_i(x)$ for $x\in S$.\\
\textbf{2.} Consider the set $T=\{x\in [n]\setminus S:h_i(x)\in f_i(S)\}$. Define the composite function $\rho=h_i(f_i^{-1}(\cdot))$. For all $x\in T$, let $h(x)=\rho^{(N)}(h_i(x))$, where $N$ is the minimum positive integer such that $\rho^{(N)}(h_i(x))\not\in f_i(S)$.\\
\textbf{3.} Let $h(x)=h_i(x)$ for all other $x\not\in S\cup T$.
\end{framed}

Regarding step \textbf{2}, it is clear that such an $N$ always exists. Indeed, observe that $|f_i(S)|$ is finite; so if there existed a $x$ such that no such $N$ existed then the sequence in $[n]$ formed by successive implementation of $\rho$ on $h_i(x)$ must be cyclic. This implies that there must exist either $x_0,x_1\in [n]$ such that $h_i(x_0)=h_i(x_1)$, or $y_0,y_1\in f_i(S)$ so that $f_i^{-1}(y_0)=f_i^{-1}(y_1)$. Both cases are absurd by virtue of the fact that both $h_i$ (on $[n]$) and $f_i^{-1}$ (on $f_i(S)$) are bijective.\\ 
 
Now, it is not too difficult to show that $\phi:h_i\mapsto h$ guarantees that $h$ is a permutation.\\

\noindent\textbf{Proof that $\phi$ gives a permutation.} Suppose $x,y\in [n]$, where $x\neq y$. Due to symmetry, there are six cases to consider.\\

\textit{Case 1.} $x,y\in S$.\\

It is obvious by step \textbf{1} that $h(x)=f_i(x)\neq f_i(y)=h(y)$.\\

\textit{Case 2.} $x\in S,y\in T$.\\

It is obvious that $h(x)\in f_i(S)$ while $h(y)\not\in f_i(S)$, so they cannot be equal.\\

\textit{Case 3.} $x\in S,y\not\in S\cup T$.\\

Same reasoning as in Case 2.\\

\textit{Case 4.} $x,y\not\in S\cup T$.\\

It is obvious by step \textbf{3} that $h(x)=h_i(x)\neq h_i(y)=h(y)$.\\

\textit{Case 5.} $x\in T,y\not\in S\cup T$.\\

Observe that $h_i(x)$, under $\rho$, gets mapped to an element in $h_i(S)$, i.e., $h(x)\in h_i(S)$. However, clearly $h(y)=h_i(y)\not\in h_i(S)$.\\  

\textit{Case 6.} $x,y\in T$.\\

Following step \textbf{2} of the procedure, let $h(x)=\rho^{(N_x)}(h_i(x))$ and $h(y)=\rho^{(N_y)}(h_i(y))$. If $N_x=N_y$ then it is straightforward that $h(x)\neq h(y)$ since $f_i^{-1}$ and $h_i$ are bijective. Now without loss of generality, suppose $N_x> N_y$, i.e., $N_x=N_y+n$ where $n\in\mathbb N$. Observe that at each iteration of $\rho$ after $N_y$ iterations of $\rho$ on $h_i(x)$ it is impossible for the sequence of integers in $[n]$ formed by each successive execution of $f_i^{-1}$ and $h_i$ at each iteration (thus it is just a length 2 sequence) to have any common members with entire sequence of integers in $[n]$ formed by the repeated execution of $f_i^{-1}$ and $h_i$, $N_y$ times, to yield $h(y)$. If not, it would imply that either $f_i^{-1}$ or $h_i$ is not bijective, a contradiction! (Drawing a diagram helps one in visualizing this.) Therefore, it is impossible for $h(x)$ and $h(y)$ to be equal in this case.\hfill $\blacksquare$\\
  
Now, we give a detailed explanation as to why $\phi$ preserves $E$ and is injective.\\

\noindent\textbf{Proof that $\phi$ preserves $E$.} By definition, $E=\bigcap_{A_{i,S}\in\mathcal{S}}\overline{A_{i',S'}}$, where $\mathcal{S}\subset\mathcal{A}\setminus\mathcal{D}_{i,S}$. Let $h_i$ satisfy $E$ and $A_{i,S}$. Under $\phi$, the positions $x\in[n]$ of the permutation $h_i$ affected are precisely only those either in $S$ or involving $h_i(S)$. Thus, if $h_i$ satisfies $\overline{A_{i',S'}}$ for each $A_{i',S'}\in\mathcal{S}$, i.e., $h_i$ does not agree with each $g_{i'}$ on each $S'$ where $S'\cap S=\emptyset$ and $g_{i}(S)\cap g_{i'}(S')$, then its image $h$ under $\phi$ also cannot agree with each $g_{i'}$ on each $S'$. Therefore, $E$ is preserved.\hfill $\blacksquare$\\

\noindent\textbf{Proof that $\phi$ is injective.} Let $h_i$ and $h'_i$ be distinct permutations satisfying both $E$ and $A_{i,S}$. This implies that they both agree with $g_i$ on $S$, i.e., $h_i(S)=h'_i(S)=g_i(S)$. For brevity, let $\phi$ map $h_i$ to $h$ and $h'_i$ to $h'$. There is some position $x_0$ at which $h_i$ and $h'_i$ differ, i.e., $h_i(x_0)\neq h'_i(x_0)$. If $x_0\in S$ or $x_0\not\in S\cup T$, we are done because $\phi$ preserves the mapping of $x_0$ via $h_i$, as well as via $h'_i$, for such $x_0$. Thus, assume $x_0\in T$. Let $h(x_0)=\rho^{(N)}(h_i(x_0))$ and $h'(x_0)=\rho^{(N')}(h'_i(x_0))$. \\

\textit{Case 1.} $N=N'$\\

Since $h_i(x_0)\neq h'_i(x_0)$, $f_i^{-1}$ maps each to distinct elements in $S$. Since $h_i$ and $h'_i$ agree on $S$ and $h_i,h'_i$ and $f_i^{-1}$ are bijective, the sequences formed by each subsequent execution of $h_i=h'_i$ and $f_i^{-1}$ respectively on $f_i^{-1}(h_i(x_0))$ and $f_i^{-1}(h'_i(x_0))$ do not share any common elements, implying that $h(x_0)\neq h'(x_0)$.\\
  
\textit{Case 2.} $N\neq N'$\\
                  
Assume, without loss of generality, that $N>N'$. Mimicking the argument presented in the first case above, one just needs to invoke the bijectivity properties of $h_i$ (on $S$) and $f_i^{-1}$ (on $f_i(S)$) to show that the two sequences cannot share any common terms even though one is longer than the other. (Drawing a diagram helps one in visualizing this.)\hfill $\blacksquare$\\

Therefore, indeed our mapping works. 
  
\end{document}